\documentclass[12pt,twoside,final,psamsfonts]{amsart}

\usepackage[psamsfonts]{amssymb}
\usepackage{times,a4wide}
\usepackage{graphicx}

\theoremstyle{plain}
\newtheorem*{theorem*}{Theorem}
\newtheorem*{TheoremA}{Theorem A}

\newtheorem{theorem}{Theorem}[section]

\newtheorem*{proposition*}{Proposition}

\newtheorem*{corollary*}{Corollary}
\newtheorem{lemma}[theorem]{Lemma}
\newtheorem*{lemma*}{Lemma}

\theoremstyle{definition}

\newtheorem*{remark*}{Remark}

\theoremstyle{definition}

\newtheorem*{definition*}{Definition}

\newcommand{\N}{{\mathbb N}}
\newcommand{\D}{{\mathbb D}}

\newcommand{\T}{{\mathbb T}}

\newcommand{\liL}{\lambda\in\Lambda}

\DeclareMathOperator{\Int}{Int}
\DeclareMathOperator{\Real}{Re}

\title[Interpolation and peak functions for the Nevanlinna and Smirnov classes]
{Interpolation and peak functions for the Nevanlinna and Smirnov classes}

\author{ Xavier Massaneda}
\address{Facultat de Matem\`atiques, Universitat de Barcelona, Catalonia}
\email{xavier.massaneda@ub.edu}

\author{ Pascal J. Thomas}
\address{Universit\'e de Toulouse, UPS, INSA, UT1, UTM, Institut de Math\'ematiques de Toulouse, F-31062 Toulouse, France}
\email{pascal.thomas@math.univ-toulouse.fr}


\thanks{The second author is supported by the Generalitat de Catalunya (grant 2009 SGR 1303) and the spanish Ministerio de Ciencia e Innovaci\'on (project MTM2011-27932-C02-01)}

\date{\today}

\begin{document}

\begin{abstract}
It is known (implicit in \cite{HMNT}) that when $\Lambda$ is an interpolating sequence for the Nevanlinna or the Smirnov class then there exist functions $f_\lambda$ in these spaces, with uniform control of their growth and attaining values 1 on $\lambda$ and 0 in all other $\lambda'\neq\lambda$. We provide an example showing that,  contrary to what happens in other algebras of holomorphic functions, the existence of such functions does not imply that $\Lambda$ is an interpolating sequence.
\end{abstract} 

\maketitle

\section{Introduction}

Consider the Nevanlinna class
\[
 N=\bigl\{f\in H(\D):\lim_{r \to 1}\int_{0}^{2\pi} \log^+|f(re^{i\theta})|\; \frac{d\theta}{2\pi}<+\infty\} ,
\]
which is a complete metric space with the distance defined by
\[
d(f,g)=\lim_{r \to 1}\int_{0}^{2\pi}  \log(1+|f(re^{i\theta})-g(re^{i\theta})|)\;  \frac{d\theta}{2\pi}\ .
\]

\begin{definition*}
A sequence $\Lambda\subset\D$ is a (free) interpolating sequence for $N$ if the space of traces $N|\Lambda$ is ideal, that is, whenever $f\in N$ and $\{\omega_\lambda\}_{\liL}$ is a bounded sequence there exists $g\in N$ such that $g(\lambda)=\omega_\lambda f(\lambda)$, $\liL$. We shall write $\Lambda\in \Int N$.
\end{definition*}

Since $N$ is an algebra, it is easily seen that $\Lambda\in \Int N$ if and only if for every bounded sequence $\{v_\lambda\}_{\liL}$ there exists $f\in N$ such that $f(\lambda)=v_\lambda $, $\liL$.
In particular, if $\Lambda\in \Int N$ there exist functions $f_\lambda\in N$ interpolating the values
\[
\delta_{\lambda,\lambda'}=
\begin{cases}
1\; & \textrm{if $\lambda^\prime=\lambda$}\\
0\;  & \textrm{if $\lambda^\prime\neq\lambda$}\ .
\end{cases}
\]
Moreover this can be achieved with functions $f_\lambda$ such that $\sup_{\liL} d(f_\lambda,0)<\infty$, as can be seen by going through the details of the proof of \cite[Theorem 1.2]{HMNT}.

Note that for other algebras of holomorphic functions the analogous size control of these $f_\lambda$ is an immediate consequence of the open mapping theorem applied to the 
restriction operator. This is the case for $H^\infty$, the algebra of bounded holomorphic functions,
the Korenblum algebra $A^{-\infty}$, or the Smirnov class $N^+$ (see below). 
However, $N$ is not even a topological vector space \cite{ShSh}, so no open mapping theorem can be used here.

Conversely, for $H^\infty$ and $A^{-\infty}$ (and for other Banach spaces), the existence of functions $f_\lambda$ with uniform control of their size and interpolating the values $\delta_{\lambda,\lambda'}$ implies that $\Lambda$ is an interpolating sequence (see \cite[Chap.VII]{Gar}, \cite[Lemma 2.3]{Ma}). We provide an example showing that this is not the case for $N$ nor $N^+$. For the Nevanlinna class, we have:

\begin{theorem}\label{thN}
Let $\Lambda$ be a sequence in $\D$.  Then, 
\begin{itemize}
\item [(a)] If $\Lambda \in \Int N$, then there exist $C>0$ and functions $f_\lambda\in N$ such that 
\begin{itemize}
\item[(i)] $d(f_\lambda,0)\leq C$ ,
\item[(ii)] $f_\lambda(\lambda^\prime)=\delta_{\lambda,\lambda'}$, $\lambda'\in\Lambda$.
\end{itemize}

\item [(b)] The converse fails: there is a sequence $\Lambda\notin \Int N$ for which there exist $C>0$ and  $f_\lambda\in N$ satisfying (i) and (ii).
\end{itemize}

\end{theorem}

On the other hand, the Smirnov class $N^+$ is defined by
\[
N^+=\bigl\{f\in N:\lim_{r \to 1}\int_{0}^{2\pi}\log^+|f(re^{i\theta})|\;  \frac{d\theta}{2\pi}=\int_{0}^{2\pi} \log^+|f^*(e^{i\theta})|\;  \frac{d\theta}{2\pi}\bigr\}\ ,
\]
where $f^*(e^{i\theta})=\lim\limits_{r \to 1} f(re^{i\theta})$ (which exists a.e. $\theta\in [0,2\pi)$).

Since $N^+$ is an $F$-space (\cite[Lemma 1]{Ya}), an application of the open mapping theorem for such spaces (see \cite[2.11, p.47]{Ru}) shows  that there exist $C>0$ and functions $f_\lambda\in N^+$ satisfying (i) and (ii) from Theorem \ref{thN}. 
But more can be said. Denote by $\mathcal F$ the class of convex, increasing functions $\psi:[0,+\infty)\longrightarrow[0,+\infty)$ such that $\lim\limits_{t\to +\infty} \psi(t)/t=+\infty$. It is known that if $f\in N^+$, there exists $\psi\in\mathcal F$,
depending on $f$, such that
\[
\int_{0}^{2\pi}
 \psi\left[  \log(1+|f^*(e^{i\theta})| \right)]\;  \frac{d\theta}{2\pi}<+\infty\ .
\]

\begin{theorem}\label{thS}
Let $\Lambda$ be a sequence in $\D$.  Then, 
\begin{itemize}
\item [(a)] If $\Lambda \in \Int N^+$, then there exist $\psi\in\mathcal F$, $C>0$ and functions $f_\lambda\in N^+$ such that 
\begin{itemize}
\item[(i)] $\displaystyle\int_{0}^{2\pi} \psi\left[  \log(1+|f_\lambda^*(e^{i\theta}) | ) \right]\;  \frac{d\theta}{2\pi}\leq C$ ,
\item[(ii)]  $f_\lambda(\lambda^\prime)=\delta_{\lambda,\lambda'}$, $\lambda'\in\Lambda$.
\end{itemize}

\item [(b)] The converse fails: there is a sequence $\Lambda\notin \Int N^+$ for which there exist  $\psi\in\mathcal F$, $C>0$ and  $f_\lambda\in N^+$ satisfying (i) and (ii).
\end{itemize}

\end{theorem}

As in the Nevanlinna case, part (a) is implicit in the proof of \cite[Theorem 1.3]{HMNT}, while part (b) will follow from an explicit example.

\section{Preliminaries. Interpolation in the Nevanlinna and Smirnov classes}

A complete description of the
interpolating sequences for $N$ and $N^+$, including a characterisation of the traces, was given in \cite [Theorems 1.2 and 1.3]{HMNT}.  In particular, we will make use of the following geometric characterisation.

Let $b_\lambda(z)=\frac{z-\lambda}{1-\bar\lambda z}$ be a Blaschke factor and 
$B_\lambda(z):= \prod_{\lambda'\neq \lambda} \left( -\frac{|\lambda'|}{\lambda'}\right) b_{\lambda'}(z)$ the Blaschke product with one factor omitted. Given a finite measure $\mu$ in $\T$, let 
\[
P[\mu](z)=\int_0^{2\pi}\frac{1-|z|^2}{|z-e^{i\theta}|^2} d\mu(\theta)\ 
\]
denote its Poisson transform. 

\begin{TheoremA} \cite{HMNT} \label{thA} 
Let $\Lambda\subset\D$. 
\begin{itemize}
 \item [(a)] $\Lambda\in \Int N$ if and only if
there exists $\mu$  positive measure with finite mass on $\T$ such that
\begin{equation*}\label{condN}
|B_\lambda(\lambda)|\geq e^{-P[\mu](\lambda)}\ ,\quad \liL .
\end{equation*}

\item [(b)]  $\Lambda\in \Int N^+$ if and only if
there exists $w\geq 0$, $w\in L^1(\T)$ , such that
\begin{equation*}\label{condS}
|B_\lambda(\lambda)|\geq e^{-P[w](\lambda)}\ ,\quad \liL .
\end{equation*}
\end{itemize}
\end{TheoremA}

In classical terminology, when  $w$ is a positive function in $L^1(\T)$, the harmonic function $u=P[w]$ is called \emph{quasi-bounded}. According to \cite[Theorem 1.3.9, p.10]{ArGa}, for any such functions there exists $\psi\in \mathcal F$ such that
\begin{equation}\label{qb}
\sup_{r<1}\int_0^{2\pi} \psi[ u (re^{i\theta})]\; \frac{d\theta}{2\pi}=\int_0^{2\pi} \psi[ w (e^{i\theta})]\; \frac{d\theta}{2\pi}< +\infty\ .
\end{equation}

\section{Proof. Necessity.}

As said before, the necessity of conditions (i) and (ii) in Theorems \ref{thN} and \ref{thS} is implicit in \cite{HMNT}. We briefly recall how this goes.

Assume $\Lambda\in \Int N$ and let $\mu$ be the measure given by Theorem A (a). Consider the function
\[
g(z)=\int_0^{2\pi}\frac{e^{i\theta}+z}{e^{i\theta}-z}\; d\mu(\theta)\ . 
\]
Since $\Real g(z)=P[\mu](z)\geq 0$ we see that $g\in N^+$ and also $ e^g\in N$.  Letting $H=(2+g)^2$ we have $H\in N^+$ and
\[
|H(\lambda)|\geq \bigl(2+\Real g(\lambda)\bigr)^2\geq \bigl(2+\log\frac 1{|B_\lambda(\lambda)|}\bigr)^2= 
 \bigl(1+\log\frac e{|B_\lambda(\lambda)|}\bigl)^2\ ,
\]
whence letting $\phi(t)=(1+t)^{-2}$ we obtain, for any $\lambda' \in \Lambda$,
\begin{align*}
|\delta_{\lambda\lambda'}| & \leq |B_\lambda(\lambda)|\ \phi \bigl(\log\frac e{|B_\lambda(\lambda)|}\bigr)\ 
e^{P[\mu] (\lambda)} |H(\lambda)|\ .
\end{align*}

\begin{theorem*} \cite[Theorem 4]{Gar2}
Let $\phi:[0,\infty)\longrightarrow [0,\infty)$ be a decreasing function such that $\int_0^\infty \phi(t)\; dt<\infty$. There exists $C>0$ such that if $\{v_\lambda\}_\lambda$ is a sequence with
\[
|v_\lambda|\leq |B_\lambda(\lambda)|\ \phi \bigl(\log\frac e{|B_\lambda(\lambda)|}\bigr) \quad \liL,
\]
then there exists $F\in H^\infty$ with $F(\lambda)=v_\lambda$, $\liL$, and $\|F\|\leq C \int_0^\infty \phi(t)\; dt$.
\end{theorem*}

Aplying Garnett's theorem to the sequence $\bigl\{ \dfrac{\delta_{\lambda\lambda'}}{e^{P[\mu] (\lambda)} |H(\lambda)| }\bigr \}_{\liL}$ we obtain a constant $C(\phi)$ and functions $F_\lambda\in H^\infty$ with $\|F_\lambda\|_\infty\leq C(\phi)$ and 
\[
F_\lambda(\lambda')=\dfrac{\delta_{\lambda\lambda'}}{e^{P[\mu] (\lambda)} |H(\lambda)| }\ .
\]
Defining
\[
f_\lambda (z)=F_\lambda(z) e^{g(z)} H(z)
\]
we finally have (ii) and
\[
\log^+|f_\lambda(z)|\leq \log C(\phi)+ P[\mu](z) +\log^+ |H(z)|\ ,
\]
which implies (i).

The same proof for the Smirnov case provides interpolating functions $f_\lambda \in N^+$ with
\[
\log^+|f_\lambda(z)|\leq \log C(\phi)+ P[w](z) +\log^+ |H(z)|\ ,
\]
where $w$ is given by Theorem A(b).
Since $H\in N^+$, the subharmonic function $\log^+ |H(z)|$ has a quasi-bouded harmonic majorant, and
 (i) in Theorem~\ref{thS} follows from \eqref{qb}.

\section{Proof. Lack of sufficiency.}
In order to construct examples of non-interpolating sequences satisfying (i) and (ii) in Theorems \ref{thN} and \ref{thS}, consider the dyadic intervals on $\T$:
\[
I_{n,k}=\{e^{2\pi i\theta} : k 2^{-n}\leq\theta< (k+1) 2^{-n}\}, \quad n\in\N,\ k=0,\dots,2^n-1\ .
\]

Let us prove first Theorem 1.1 (b).

Consider the sequence $A$ defined in the following way: on the ray terminating at an end of a dyadic interval of the $n$-th generation  (and which is not and end of an interval in a previous generation) consider the dyadic sequence with radii $1-2^{-m}$, starting at $m=2n$. Explicitly
\[
A=\{a_{m}^{n,k}\}_{
\begin{subarray}{l}
n\in\N\\
0\leq k<2^n-1\textrm{, $k$ odd}\\
m\geq 2n
\end{subarray} }
\qquad a_{m}^{n,k}=(1-2^{-m}) e^{2\pi i k2^{-n}}\ .
\]

\begin{figure}[htb]
\centering
\includegraphics[width=0.95\textwidth]{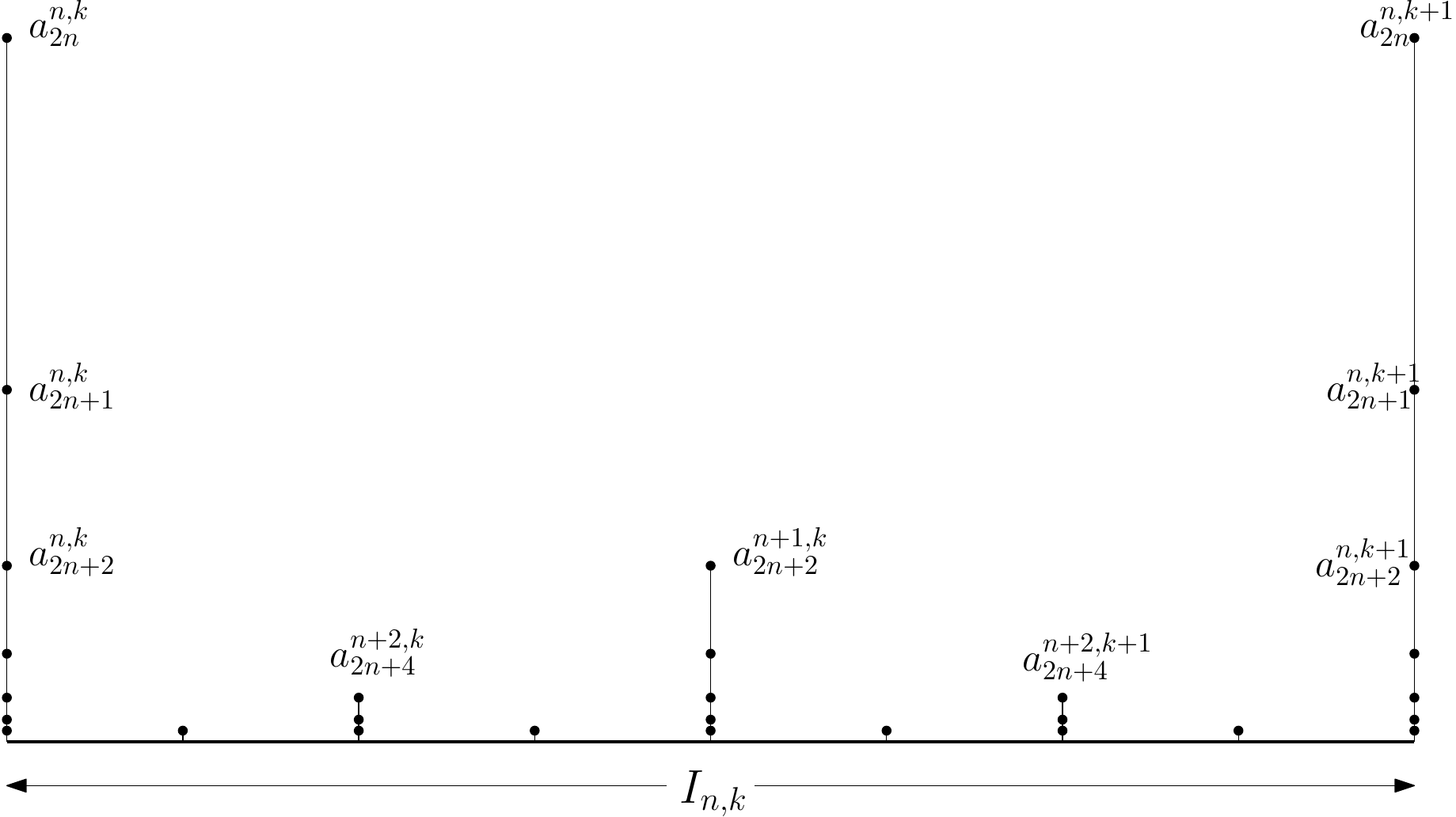}
\caption{Representation of the points of $A$ ``above'' the interval $I_{n,k}$}
\label{fig:awesome_image}
\end{figure}

Now, to each $a\in A$ associate a point $b$ on the same ray and so that 
\begin{equation}\label{sepN}
\varrho(a,b):=\left|\frac{a-b}{1-\bar a b}\right|\simeq \exp\left(-\frac{1}{1-|a|}\right)\ .
\end{equation}
We thus obtain a sequence $B$, which can be explicitly given by
\[
B=\{b_{m}^{n,k}\}_{
\begin{subarray}{l}
n\in\N\\
 0\leq k<2^n-1 \textrm{, $k$ odd}\\
m\geq 2n
\end{subarray} }
\qquad b_{m}^{n,k}=(1-e^{-2^m})a_{m}^{n,k}\ .
\]

\begin{lemma}\label{intH}
The sequences  $A$ and $B$ are both interpolating for $H^\infty$.
\end{lemma}

Recall that $A\in \Int H^\infty$ means that there exists $C>0$ (the interpolation constant) such that for every bounded sequence $\{v_a\}_{a\in A}$ there is $F\in H^\infty$ with $\|F\|_\infty\leq C \|\{v_a\}\|_\infty$ and $F(a)=v_a$, $a\in A$.

\begin{proof}
Let us see that $A\in \Int H^\infty$.
By Carleson's theorem (\cite[Theorem 1.1, Chap.VII]{Gar}) it is enough to see that $A$ is separated in the pseudo-hyperbolic metric $\varrho$, which is obvious from the definition, and that $\nu=\sum_{a\in A} (1-|a|) \delta_a$ is a Carleson measure.

Let us see first that $A$ is a Blaschke sequence:
\begin{align*}
\sum_{a\in A} 1-|a|=\sum_{n=1}^\infty \sum_{k=0}^{2^n-1}\sum_{m\geq 2n} 2^{-m}
\simeq \sum_{n=1}^\infty 2^n 2^{-2n}<+\infty
\end{align*}

In order to see that $\nu$ is a Carleson measure we have to prove that there exists $C>0$ such that $\nu(Q(I))\leq C |I|$ for all $I$ interval in $\T$, where $Q(I)=\{re^{i\theta} : 1-r<|I|\ ,\, e^{i\theta}\in I\}$ is the associated Carleson box.  It is enough to consider the case where $I$ a dyadic interval. Thus let $I=I_{n,k}$ and
\[
Q(I_{n,k})=\{re^{i\theta} : r>1-2^{-n}\ ,\, \theta\in 2\pi[k2^{-n}, (k+1)2^{-n})\}\ .
\]
By construction, in $Q(I_{n,k})$ there are  $2^j$ rays of the $(n+j)$-th generation. Hence
\[
\sum_{a\in A\cap Q(I_{n,k})} 1-|a|\simeq \sum_{j=1}^\infty 2^j \sum_{m\geq 2(n+j)} 2^{-m}\simeq
\sum_{j=1}^\infty 2^{-2n-j}\lesssim 2^{-n} =|I_{n,k}|\ .
\]
\end{proof}

Define $\Lambda=A\cup B$.
Let us see first that there exist $C>0$ and $f_\lambda$ satisfying (i) and (ii) in Theorem~\ref{thN}.

Fix $\lambda=a_m^{n,k}$ and denote by $\tilde\lambda=b_m^{n,k}$ its ``twin". As just seen, there exist $C>0$ and $P_\lambda^A\in H^\infty$ such that
\begin{itemize}
\item $\| P_\lambda^A \|_\infty \leq C$,
\item $P_\lambda^A (\lambda)=1$ and $P_\lambda^A(a_{m^\prime}^{n^\prime,k^\prime}) =0\quad \forall(n^\prime,k^\prime,m^\prime)\neq (n,k,m)$\ .
\end{itemize}

As in the proof of Lemma~\ref{intH}, we can see that $(B\setminus\{\tilde\lambda\})\cup \{\lambda\}$ is also in $\Int H^\infty$, and with interpolation constant $C>0$ independent of $\lambda$. Therefore there exist  $P_\lambda^B\in H^\infty$ such that
\begin{itemize}
\item $\| P_\lambda^B \|_\infty \leq C$,
\item $P_\lambda^B (\lambda)=1$ and $P_\lambda^B(b_{m^\prime}^{n^\prime,k^\prime}) =0\quad \forall(n^\prime,k^\prime,m^\prime)\neq (n,k,m)$\ .
\end{itemize}

Define finally
\begin{equation}\label{fl}
f_\lambda:= c_\lambda P_\lambda^A  P_\lambda^B b_{\tilde\lambda} e^{g_\lambda}\ ,
\end{equation}
where 
\[
g_\lambda(z)=\frac{\lambda^*+z}{\lambda^*-z} \qquad  (\lambda^*=\lambda/|\lambda|)\ 
\]
and $c_\lambda$ is chosen so that $f_\lambda(\lambda)=1$.

Notice that, by construction, (ii) in Theorem~\ref{thN} holds. In order to see (i) notice that \eqref{sepN} gives
\[
|c_\lambda|=\frac{1}{|b_{\tilde\lambda}(\lambda)|}\exp\left(-\frac{1+|\lambda|}{1-|\lambda|}\right)\simeq
\exp\left(\frac 1{1-|\lambda|}- \frac{1+|\lambda|}{1-|\lambda|}\right)\lesssim 1\ .
\]
Then
\begin{align*}
\log|f_\lambda(z)|=
\log|c_\lambda|+\log|P_\lambda^A(z)|+\log|P_\lambda^B(z)|+\log| b_{\tilde\lambda}(z)|+\Real g_\lambda(z)
\end{align*}
and therefore
\[
\log^+|f_\lambda(z)|\leq \log \|c\|_\infty + 2\log C+ \Real g_\lambda(z)\ .
\]
Since
\[
\sup_{r<1}\int_0^{2\pi}   \Real g_\lambda(r e^{i\theta}) \frac{d\theta}{2\pi}=
\sup_{r<1} \int_0^{2\pi} \frac{1-r^2}{|r e^{i\theta}-\lambda^*|^2}\frac{d\theta}{2\pi}=1\  ,
\]
we have (i), as desired.

Let us see now that $\Lambda\notin \Int N$ by seeing that there is no $\mu$ satisfying the condition of Theorem A(a). Since
\[
\log\frac 1{|B_\lambda(\lambda)|}\simeq \frac 1{1-|\lambda|}\qquad \liL ,
\]
such $\mu$ should satisfy in particular  (fixed any $n,k$)
\[
1\lesssim (1-|a_{m}^{n,k}|)  P[\mu] (a_{m}^{n,k}), \quad \forall m\geq 2n\ .
\]
This would force the measure $\mu$ to satisty $1\lesssim \mu\{e^{2\pi ik2^{-n}}\} $ for all $n,k$ (\cite[Theorem 2.2]{ShSh}), and therefore $\mu(\T)$ could not be finite.

Let us prove now Theorem 1.2 (b). In the same construction done for $N$ consider a sequence made of the ``first" couple of points of each ray, and with a slightly bigger separation. More precisely, let $\tilde\Lambda=\tilde A\cup\tilde B$, where
\[
\tilde A=\{a_{n,k}\}_{
\begin{subarray}{l}
n\in\N\\
0\leq k<2^n-1,\ \textrm{$k$ odd}
\end{subarray} }
\qquad a_{n,k}=(1-2^{-2n}) e^{2\pi i k2^{-n}}\ .
\]
and $\tilde B=\{b_{n,k}\}_{n,k}$ is so that $a_{n,k}$ and $b_{n,k}$ are on the same ray and 
\begin{equation}\label{sepS}
\varrho(a_{n,k},b_{n,k})=\exp\left(-\frac 1{(1-|a_{n,k}|)\log_2(\frac 1{1-|a_{n,k}|})}\right) =\frac{2^{2n}}{2n}\ .
\end{equation}

As in Lemma~\ref{intH}, $\tilde A$ and $\tilde B$ are $H^\infty$-interpolating sequences, so there exist bounded peak functions $P_\lambda^{\tilde A}$, $P_\lambda^{\tilde B}$ with the same properties as before. 

Given $\lambda\in\tilde\Lambda$, let $I_\lambda$ denote the Privalov ``shadow" of $\lambda$ on $\T$, that is
\[
I_\lambda=\{e^{i\theta} : |e^{i\theta}-\lambda|\leq 2(1-|\lambda|)\}\ .
\]
Let 
\[
w_\lambda(\theta)=\frac {C_0}{(1-|\lambda|)\log_2(\frac 1{1-|\lambda|})}  \chi_{I_\lambda}(e^{i\theta}) , 
\]
where $c$ is a universal constant to be chosen later, and consider
\[
g_\lambda(z)=\int_0^{2\pi} \frac{e^{i\theta}+z}{e^{i\theta}-z}\; w_\lambda(\theta)\; \frac{d\theta}{2\pi}\ .
\]
Notice that $\|w_\lambda\|_{L^1(\T)}\simeq c (\log_2(\frac 1{1-|\lambda|}))^{-1}\lesssim 1$ and
\[
\Real g_\lambda(\lambda)=P[w_\lambda](\lambda)\simeq\frac 1{1-|\lambda|}\int_{I_\lambda} w_\lambda(\theta)\; d\theta \ge \frac {C_1 C_0}{(1-|\lambda|)\log_2(\frac 1{1-|\lambda|})}\ .
\]

Now define $f_\lambda$ as in \eqref{fl}, with these new $g_\lambda$. Again, it is clear that (ii) in Theorem~\ref{thS} holds. 
Also, $\{c_\lambda\}_{\liL}$ is bounded if $C_0$ is chosen appropriately: if $\lambda=a_m^{n,k}$,
\[
\log|c_\lambda|=\log\frac 1{\varrho(a_{n,k},b_{n,k})}-\Real g_{a_{n,k}}(a_{n,k})\leq\frac {1-{C_1 C_0}}{(1-|a_{n,k}|)\log_2(\frac 1{1-|a_{n,k}|})}\leq 0\ .
\]
In order to see (i) notice that, as before, there exists $\tilde C>0$ such that
\begin{align*}
\log^+|f_\lambda|&\leq \log \|c\|_\infty+2 \log C +P[w_\lambda]\leq \tilde  C +P[w_\lambda]\ .
\end{align*}
Therefore, by \eqref{qb}, for $\lambda=a_m^{n,k}$ and taking $\psi(t)=(1+t)\log(1+t)\in\mathcal F$ we get (i):
\begin{align*}
\int_0^{2\pi} \psi[\log^+|f_\lambda^*(e^{i\theta})|]\, \frac{d\theta}{2\pi}&\leq 
\int_0^{2\pi} \psi[ \tilde C+\frac {{C_0}\chi_{I_\lambda}(\theta)}{(1-|\lambda|)\log_2(\frac 1{1-|\lambda|})} ]\, \frac{d\theta}{2\pi}\\
&=\int_{\theta\notin I_\lambda} \psi(\tilde  C) \, \frac{d\theta}{2\pi}
+\int_{I_\lambda}\psi \bigl[\tilde C+\frac {C_0}{(1-|\lambda|)\log_2(\frac 1{1-|\lambda|})}\bigr] \, \frac{d\theta}{2\pi}\\
&\lesssim 1+ (1-|\lambda|) \psi \bigl[\frac {2{C_0}}{(1-|\lambda|)\log_2(\frac 1{1-|\lambda|})}\bigr]\\
&\lesssim 1+\frac {1}{\log_2(\frac 1{1-|\lambda|})}\log_2 \left(\frac {2{C_0}}{(1-|\lambda|)\log_2(\frac 1{1-|\lambda|})}\right)
\lesssim 1\ .
\end{align*}

Let us finish by proving that $\Lambda\notin \Int N^+$. Assume that there is $w\in L^1(\T)$ satisfying Theorem A (b). Then
\begin{align*}
\frac 1{(1-|a_{n,k}|)\log_2(\frac 1{1-|a_{n,k}|})}&=\log\frac 1{\varrho(a_{n,k}, b_{n,k})}\leq \log\frac 1{|B_{a_{n,k}}(a_{n,k})|}\\
&\leq P[w](a_{n,k})=\int_0^{2\pi}\frac{1-|a_{n,k}|^2}{|a_{n,k}-e^{i\theta}|^2}\; w(\theta)\; \frac{d\theta}{2\pi}\ ,
\end{align*}
and therefore
\[
\sum_{n\geq 1} \sum_{k=0}^{2^n-1}\frac 1{2n}\simeq \sum_{n\geq 1} \sum_{k=0}^{2^n-1}\frac 1{\log_2(\frac 1{1-|a_{n,k}|})}\lesssim
\int_0^{2\pi}\sum_{n,k}\frac{(1-|a_{n,k}|^2)^2}{|a_{n,k}-e^{i\theta}|^2}\; w(\theta)\; \frac{d\theta}{2\pi}\ .
\]
We will have a contradiction as soon as we prove that
\[
\sup\limits_{\theta\in[0,2\pi)}\sum_{n,k}\frac{(1-|a_{n,k}|^2)^2}{|a_{n,k}-e^{i\theta}|^2}<\infty\ .
\]
With no loss of generality assume that $e^{i\theta}=1$ and that the dyadic intervals $I_{n,k}$, $0\leq k<2^n-1$,  are ordered so that $e^{i\theta}\in I_{n,0}$. Then we have 
\begin{align*}
|e^{i\theta}-a_{n,k}|^2 &=|1-a_{n,k}|^2\simeq (1-|a_{n,k}|)^2 + |e^{2\pi i k2^{-n}}-1|^2\\
& \simeq (2^{-2n})^2 + (k 2^{-n})^2= 2^{-2n} (2^{-2n}+k^2)\ .
\end{align*}
Since for each $e^{i\theta}$ there is only one $a_{n,k}$ with $1-|a_{n,k}|\simeq |e^{i\theta}-a_{n,k}|$, we have then
\begin{align*}
\sum_{n,k}\frac{(1-|a_{n,k}|^2)^2}{|a_{n,k}-e^{i\theta}|^2}&\simeq
\sum_{n\geq 1} \sum_{k=1}^{2^n-1} \frac{2^{-4n}}{2^{-2n} k^2}\simeq \sum_{n\geq 1} 2^{-2n}\ .
\end{align*}


\begin{thebibliography}{RosROV}

\bibitem[ArGa]{ArGa}
David H. Armitage, Stephen J. Gardiner, {Classical potential theory}, Springer Monographs in Mathematics (2000).

\bibitem[Gar77]{Gar2}
John B. Garnett, ``Two remarks on interpolation by bounded analytic functions'', Banach spaces of analytic functions (Proc. Pelczynski Conf., Kent State Univ., Kent, Ohio, 1976), pp. 32--40. Lecture Notes in Math., Vol. \textbf{604}, Springer, Berlin (1977). 


\bibitem[Gar07]{Gar}
John B. Garnett, {Bounded analytic functions. Revised first edition}, Graduate Texts in Mathematics, 236. Springer, New York, 2007.

\bibitem[HMNT]{HMNT}
Andreas Hartmann, Xavier Massaneda, Artur Nicolau, and Pascal J. Thomas. 
``Interpolation in the Nevanlinna class and harmonic majorants'', J. Funct. Anal. \textbf{217} (2004) 1--37. 

\bibitem[Ma]{Ma}
Xavier Massaneda. 
``Interpolation by holomorphic functions in the unit ball with polynomial growth'', Ann. Fac. Sci. Toulouse Math. \textbf{6} (6) (1997) 277--296. 

\bibitem[Ru73]{Ru}
Walter Rudin, {Functional Analysis}, McGraw-Hill, New York, 1973.


\bibitem[ShSh]{ShSh}
Joel S.  Shapiro, Alan L. Shields,  ``Unusual topological properties of the Nevanlinna class'', Amer. J. of Math., \textbf{97}, n.4 (1975)  915--936.

\bibitem[Ya]{Ya}
Niro Yanagihara,  ``Interpolation theorems for the class $N^+$'', Illinois J. Math.,\textbf{18},  (1974)  427--435.

\end{thebibliography}
\end{document}